\documentclass[12pt]{amsart}
\usepackage{times}
\usepackage[margin=1.2in]{geometry}
\usepackage[T1]{fontenc}
\usepackage[utf8]{inputenc}
\usepackage{microtype}
\usepackage{hyperref}

\usepackage{amsmath,amssymb,amsthm,mathtools}
\usepackage{tikz-cd}

\theoremstyle{plain}
\newtheorem{theorem}{Theorem}[section]
\newtheorem{proposition}[theorem]{Proposition}
\newtheorem{lemma}[theorem]{Lemma}
\newtheorem{corollary}[theorem]{Corollary}

\theoremstyle{definition}
\newtheorem{definition}[theorem]{Definition}
\newtheorem{example}[theorem]{Example}

\theoremstyle{remark}
\newtheorem{remark}[theorem]{Remark}

\newcommand{\N}{\mathbf{N}}
\newcommand{\Z}{\mathbf{Z}}

\newcommand{\F}{\mathbf{F}}
\newcommand{\w}{\mathfrak{w}}
\renewcommand{\epsilon}{\varepsilon}
\newcommand{\act}{\curvearrowright}
\newcommand{\deriv}{\mathrm{d}}

\DeclareMathOperator{\im}{im}
\DeclareMathOperator{\outd}{outdeg}
\DeclareMathOperator{\ind}{indeg}
\DeclareMathOperator{\dom}{dom}

\usepackage{xcolor}

\title{Elasticity of free type III actions of free groups}
\author{Antoine Poulin}
\date{}

\usepackage[english]{babel}

\begin{document}

\begin{abstract}
    We prove that measure-class-preserving non-amenable treeable equivalence relations of type III, meaning not preserving any equivalent $\sigma$-finite measure, are induced by free actions of non-abelian free groups of any given number of generators, including infinitely generated free groups, with the additional property that no ends of the induced Schreier graph are vanishing. This is done using a characterization of type III due to Hopf for transformations and Dang-Ngoc-Nghiem in general. This highlights the difference between the type III setting and the measure-preserving setting.
\end{abstract}

\maketitle
\section{Introduction}

Measured group theory is the study of groups through their actions on measure spaces. Much information about a probability measure-preserving (\textbf{pmp}) action is recoverable solely its orbit equivalence relation. Two pmp actions $\Gamma \act (X,\mu)$ and $\Delta \act (Y,\nu)$ of countable groups on standard probability spaces are said to be \textbf{orbit equivalent (OE)} if there is a measure isomorphism from $(X,\mu)$ to $(Y,\nu)$ preserving the orbit equivalence relations. A result of Dye (\cite{dye_groups_1959}) says that any two ergodic pmp $\Z$-actions are OE, and a subsequent theorem of Ornstein and Weiss (\cite{ornstein_ergodic_1980}) states that every ergodic pmp action of an amenable group is OE to an ergodic pmp $\Z$-action. 

Cost is an OE-invariant of pmp equivalence relations defined as half the infimum over graphings of their expected degrees, where a graphing of an equivalence relation is a Borel graph whose connected components are exactly the equivalence classes. It was introduced by Levitt in \cite{levitt_cost_1995} and studied extensively by Gaboriau in \cite{gaboriau_mercuriale_1998, gaboriau_cout_2000} and subsequent works. A question of Levitt in  \cite{levitt_cost_1995} was whether cost is a non-trivial OE-invariant for free actions of free groups. Gaboriau answered this in \cite{gaboriau_mercuriale_1998, gaboriau_cout_2000}, in fact proving the stronger result that free groups of differing ranks are not OE. The main result of \cite{gaboriau_mercuriale_1998, gaboriau_cout_2000} is that cost is attained by a \textbf{treeing}, an acyclic Borel graphing on $X$, if it exists. This implies that the cost of any free pmp $\F_d$-action is exactly $d$. A converse was also given by Hjorth in \cite{hjorth_lemma_2006}. Using edge sliding technique, he shows that 
a graphing of cost greater than $d$ can be regularized so that it contains $d$ transformations. Combined with the work of Gaboriau in \cite{gaboriau_cout_2000}, one gets the corollary that if an ergodic pmp relation admits a treeing and has cost $d$, there is an ergodic pmp $\F_d$-action OE to it. 

It directly follows from \cite[Theorem 3.17]{JKL} that given a compressible treeable countable Borel equivalence relation, there is a free action of $\F_d$ generating this equivalence relation, for any $d\in \{2, ..., \infty\}$. By a theorem of Nadkarni \cite{Nadkarni}, compressibility is equivalent to the inexistence of
invariant probability measures.

A recent program is to transfer measurable graph combinatorics results from the pmp setting to the measure-class-preserving (\textbf{mcp}) setting, where the equivalence relation may not preserve the probability measure, but preserves null sets. This was initiated in \cite{Tserunyan_2022}, where Tserunyan constructs hyperfinite ergodic subgraphs of ergodic mcp graphs. Since then, many results have been generalized to the mcp setting. For example, in \cite{bowen2022oneended}, Bowen, Zomback and the author show the existence of one-ended spanning treeings in measure-class-preserving actions of one-ended amenable group. A characterization of amenability for mcp acyclic graphs has recently been given by Tserunyan and Tucker-Drob in \cite{TTD}, generalizing the characterization of Adams in \cite{Adams_1990}. This states that an mcp acyclic graph is amenable if and only if a.e component has at most 2 non-vanishing ends. This has then been successfully applied in \cite{chen2023nonamenablesubforestsmultiendedquasipmp} by Chen, Terlov and Tserunyan to prove the existence of non-amenable subforests of locally finite graphs with infinitely many non-vanishing ends, akin to a result due to Gaboriau and Ghys (See \cite{Ghys} and \cite[IV.24]{gaboriau_cout_2000}) in the pmp setting. 

In light of \cite{TTD}, it is clear that non-vanishing of ends is a non-triviality condition. Even though it follows from \cite{JKL} that compressible relations can be generated by free actions of non-abelian free groups of any rank, it is a priori unclear whether the ends of such actions can be made non-vanishing when equipped with an invariant measure class. Our main result is the following ergodic-theoretic strengthening in the \textbf{type III} setting,  where the measure class does not contain an invariant ($\sigma$-finite) measure.

\begin{theorem}
        Let $(X, \mu, G)$ be an ergodic type III non-amenable locally countable Borel acyclic graph and $d \in \{2, \dots, \infty\}$. Then there is a free Borel action of $\F_d$ on $X$ such that $E_G = E_{\F_d}$ a.e. and a.e the components of the Schreier graph does not have vanishing ends.
\end{theorem}

This is Theorem \ref{theo:main} in the text. This answers a question of D. Gaboriau, A. Tserunyan and R. Tucker-Drob. The techniques are inspired by Hjorth's edge sliding argument, further developed in \cite{miller2018edgeslidingergodichyperfinite}. Conversely, see Example \ref{examp:IIinf} for a compressible graph with an invariant infinite measure which does not admit $\F_\infty$ graphing without vanishing ends for any equivalent probability measures.

Informally, this result also suggests that if a form of cost exists for measure-class-preserving actions, it will be very different than in the pmp setting (see also Remark \ref{rem:triv}). These results exhibit the fundamental differences between the type III and the pmp setting.

Despite the seemingly negative overtone of our result to the tractability of the type III setting, in upcoming work with D. Gaboriau, A. Tserunyan, R. Tucker-Drob and K. Wr\'{o}bel on the orbit equivalence classification of Baumslag--Solitar groups, we crucially use the type III setting to solve a problem in the pmp setting, in particular using the results proven here.



\subsection*{Acknowledgements}

The author would like to thank Damien Gaboriau, Alexander Kechris, Sam Mellick, Anush Tserunyan, Robin Tucker-Drob and Konrad Wróbel for helpful discussions. Special thanks go to Damien Gaboriau and Anush Tserunyan for their guidance and support.

\section{Preliminaries on measured combinatorics}

In this section, we establish some preliminaries on measured equivalence relations and Borel graphs. The material is standard, with references such as \cite{kechris_topics_2004} and \cite{Tserunyan_2022}.

A \textbf{standard measure space} $(X,\mu)$ is a Polish space $X$ equipped with a $\sigma$-finite Borel measure $\mu$ on $X$. All measures we consider will be $\sigma$-finite and Borel, i.e supported on the Borel $\sigma$-algebra. Equality of sets and functions is always considered mod null. Two $\sigma$-finite Borel measures on $X$ are \textbf{equivalent}, denoted $\mu' \sim \mu$, if for all Borel subset $A \subset X$,
\[ \mu'(A) = 0 \Longleftrightarrow \mu(A) =0 .\]
In this case, there exists an essentially unique Radon-Nikodym derivative $\frac{\deriv \mu'}{\deriv \mu}: X \rightarrow (0, \infty)$ satisfying 
\[\int_X f(x) \; \deriv \mu'(x) = \int_X f(x) \cdot \frac{\deriv \mu'}{\deriv \mu}(x) \; \deriv \mu(x),\]
for all Borel function $f: X \rightarrow [0,\infty]$.

A \textbf{countable Borel equivalence relation} on a standard measure space $(X,\mu)$ is an equivalence relation $E \subset X^2$ which is Borel as a subset of $X^2$ and whose equivalence classes are all countable. We will always abbreviate ``countable Borel equivalence relation'' to \textbf{CBER}. The $E$-\textbf{saturation} of a subset $A \subset X$ is
\[ [A]_E = \left\{x \in X : \exists y \in A, x \mathbin{E} y \right\}.\]
We write $[x]_E$ for $[\{x\}]_E$, the equivalence class of $x$. A subset $A \subset X$ is a \textbf{complete section} if $[A]_E$ is co-null. We say $E$ is \textbf{ergodic} if whenever $A \subset X$ is non-null, $A$ is a complete section. We further say that it is \textbf{properly ergodic} if $\mu$ is not supported on a single $E$-class. We define two measures on $E$: for $A \subset E$, let
\begin{align*}
    \nu_l(A) &= \int_X \left|A \cap s^{-1} \{x \}\right|\; \deriv \mu\\
    \nu_r(A) &= \int_X \left|A \cap t^{-1} \{x \}\right|\; \deriv \mu
\end{align*}
where $s, t: X^2 \rightarrow X$ denote the two projections. We say that $E$ is \textbf{measure-preserving} if $\nu_l = \nu_r$ and say that $E$ is \textbf{measure-class-preserving} if $\nu_l \sim \nu_r$. We abbreviate these as \textbf{mp} and \textbf{mcp}, respectively.

\begin{remark}
    In an mcp CBER $(X,\mu, E)$, the measure class of $\nu_l$ and $\nu_r$ allows us to speak of non-null sets of edges. For ergodic CBERs, proper ergodicity is equivalent to this measure class being atomless.
\end{remark}

\begin{definition}
    Let $(X,\mu, E)$ be an mcp CBER. The \textbf{Radon-Nikodym cocycle} is 
    \[ \w(y,x) := \frac{\deriv \nu_l}{\deriv \nu_r}(y,x).\]
    It satisfies the cocycle identity
    \[ \w(z,y) \cdot \w(y,x) = \w(z,x).\]
\end{definition}

Informally, we think of this cocycle as giving a relative weight between two points of the same class. In particular, a measure is invariant if and only if its Radon-Nikodym cocycle is equivalently $1$. This is formalized through the mass-transport principle \cite[Lemma 5.2]{Tserunyan_2022}.

For an mcp CBER $(X, \mu, E)$, the \textbf{full pseudogroup} associated to $E$, denoted $[[E]]$, is the set of Borel bijections $\phi:A \rightarrow B$ with Borel domain and image satisfying $\phi(x) \mathrel{E} x$ for all $x \in X$. We identify elements of $[[E]]$ if they agree mod null.

If $G$ is any locally countable Borel graph on $X$, denote $E_G$ for the CBER induced by $G$. For a locally countable Borel graph $G$ on $X$ and a point $x\in X$, we write $[x]_G$ for $[x]_{E_G}$. Similarly, we say that $G$ is ergodic if $E_G$ is, and similarly for other notions. If $H \leq G$ is a subgraph, we say it is a \textbf{strict} subgraph and write $H<G$ if $E_H \neq E_G$ on a positive measure set.

An \textbf{end} of the connected component of an acyclic graph is a class of geodesics $(x_0, x_1, ...)$ mod tail equivalence. Given an acyclic, locally countable Borel graph $G$ on $X$, and $(x,y) \in G$ a directed edge, denote $N_{(x,y)} = \{z \in [x]_{G} : d_G(y,z) < d_G(x,z)\}$, where $d_G$ is the graph distance. Say that an end represented by a geodesic $(x_0, x_1, ...)$ is \textbf{vanishing} if 
\[ \limsup_{y\in N_{(x_{i},x_{i+1})}} \w(y, x_0 ) \longrightarrow 0.\]

 It is straightforward to check that an end vanishing does not depend on the chosen representative. See \cite{TTD} for a more complete investigation of vanishing and non-vanishing ends. Of note, an acyclic locally countable Borel graph on a standard probability space is smooth, i.e admits a Borel selection of one vertex per components on a conull set, if and only if every end of a.e every components vanishes.
 
 The fact that every standard Borel space admits a Borel linear ordering implies:

\begin{lemma}\label{lem:order}
    Let $(X, \mu)$ be a standard measure space and $G$ a Borel graph on $X$. Then, there is a Borel orientation $\overrightarrow{G} \subset G$.
\end{lemma}

Here, an orientation means a unique, coherent choice of $(x,y)$ or $(y,x)$ for every $(x,y) \in G$.

\subsection{Type of a measured equivalence relation}

In this section, we recall the classification of mcp CBERs in types and establish some relevant properties. The material here is less standard, but still folklore. It is mainly based on work of W. Krieger, for example \cite{krieger_non-singular_1969,krieger_non-singular_1969-1}. See \cite{hamachi1981ergodic} and \cite{Katznelson_Weiss_1991} for references.

\begin{definition}
An ergodic mcp CBER $(X, \mu, E)$ is:
\begin{itemize}
		\item \textbf{type I} if $\mu$ is atomic. Equivalently, if it is not properly ergodic.
    \item \textbf{type II} if there is a $\sigma$-finite $\mu' \sim \mu$ preserved by $E$ and $\mu$ is non-atomic.
    \item \textbf{type III} if there is no $\sigma$-finite $\mu' \sim \mu$ preserved by $E$.
\end{itemize}
\end{definition}

\begin{proposition}\label{prop:cohom}
    A properly ergodic mcp CBER $(X,\mu, E)$ is type II if and only if its Radon-Nikodym cocycle is a coboundary, i.e there is a Borel function $w:X \rightarrow (0,\infty)$ for which
    \[  w(y) = \w(y,x) \cdot  w(x).\]
\end{proposition}
\begin{proof}[Sketch of proof]
    Notice that when $\mu' \sim \mu$, the Radon-Nikodym cocycles satisfy
    \[\w_{\mu'}(y,x) = \frac{\deriv \mu}{\deriv \mu'}(x) \cdot \w_\mu(y,x) \cdot \frac{\deriv \mu'}{\deriv \mu}(y) .\]
    Consider the function $w$ as the Radon-Nikodym derivative with respect to some preserved $\mu' \sim \mu$, i.e $w(x) = \frac{\deriv \mu}{\deriv \mu'}(x)$. 
\end{proof}

In terms of the informal interpretation of the cocycle as a relative weight function, Proposition \ref{prop:cohom} states that a CBER is type II if its Radon-Nikodym cocycle is derived from an absolute weight function.

\begin{corollary}\label{cor:stable}
    Type II and type III are preserved under taking complete sections: If $(X, \mu, E)$ is a properly ergodic mcp type II (resp. type III) CBER and $A \subset X$, then the restriction $(A, E|_A, \mu|_A)$ is also type II (resp. type III).
\end{corollary}
\begin{proof}[Sketch of proof]
    Restrict or extend $w$ as in Proposition \ref{prop:cohom}. To extend, notice that if $x \not\in A$, there is a unique $w_x$ such that $w_x \cdot \w(a,x) = w(a)$, for all $a \in [x]_E \cap A$.
\end{proof}

If $(X, \mu, E)$ is an mcp CBER, $x \in X$ and $B \subset [x]_E$, write $\w(B,x) = \left\{\w(b,x) : b \in B \right\}$. 

\begin{corollary}\label{cor:unbound}
     Let $(X,\mu, E)$ be a properly ergodic CBER. Then $E$ is type III if and only if for all $A\subset X$ non-null and for all $x \in X$, the values of the cocycle $\w( [x]_E \cap A, x)$ are unbounded above and below (in $(0,\infty)$). 
\end{corollary}
\begin{proof}
    ($\Longrightarrow$) Suppose first $E$ is type III. By Corollary \ref{cor:stable}, it suffices to show that the cocycle is unbounded on $A = X$. Suppose towards a contradiction that for all $x \in X$,
    \[\frac{1}{w(x)} := \inf\left\{\w(z,x) : z \in [x]_E \right\} > 0.\]
    Notice that if $y \mathbin{E} x,$
    \begin{align*}
        \frac{1}{w(x)}  &= \inf\left\{\w(z,x) : z \in [x]_E \right\} \\
             &= \inf\left\{\w(z,y) \cdot \w(y,x): z \in [x]_E \right\} \\
             &= \w(y,x) \cdot \inf\left\{\w(z,y) : z \in [x]_E \right\} \\
             &= \w(y,x) \cdot  \frac{1}{w(y)} ,
    \end{align*}
    as required for Proposition \ref{prop:cohom}, contradicting that $E$ type III. Similarly, use $\sup$ if the cocycle is bounded above.

    ($\Longleftarrow$) Suppose now that $E$ is type II and let $w$ as in Proposition \ref{prop:cohom}. Let $R$ be large enough so that 
    \[ A = \{x \in X: w(x) \in (R^{-1}, R) \}\]
    is non-null. Then, for $x,y \in A$ with $x \mathbin{E} y$,
    \begin{align*}
        \w(y,x) = \frac{w(x)}{w(y)} \in (R^{-2}, R^2).
    \end{align*}
    Thus, on $A$ the cocycle is bounded above and below.
    \end{proof}

The following characterisation of type III is the key tool in our result. It follows from \cite{hopf_theory_1932} for hyperfinite CBER and \cite{dang-ngoc-nghiem_classification_1973} for general CBER.

\begin{proposition}\label{prop:bij-set}
Let $(X,\mu, E)$ be a properly ergodic type III CBER. Let $A, B \subset X$ Borel non-null sets. Then there exists $\phi \in [[E]]$ with $\dom(\phi) = A$, $\im(\phi) = B$.
\end{proposition}

\section{Sliding edges along a type III subgraph}

\begin{lemma}\label{lem:bij-edge}
    Let $(X,\mu, E)$ be an ergodic type III CBER and $F \subset E$ an ergodic type III subrelation. Let $\mathcal{A} \subset E - \Delta$ be a Borel non-null collection of directed edges and $A \subset X$ a Borel non-null subset of vertices. Then there is a Borel bijection $\phi:\mathcal{A} \rightarrow A$ such that $\phi(e) \mathbin{F} s(e)$ for all $e \in \mathcal{A}$.
\end{lemma}
\begin{proof}
    Consider the CBER $(E, \nu_l, F^{(2)})$ where
    \[ e \mathbin{F^{(2)}} e' \longleftrightarrow s(e) \mathbin{F} s(e').\]
    Embedding $X$ as the diagonal $\Delta$ gives rise to an isomorphism $(\Delta, \nu_l, F^{(2)}) \cong (X, \mu, F)$. Since $(X, \mu, F)$ is type III, so is $(E, \nu_l, F^{(2)})$, by Corollary \ref{cor:stable}.

    Apply Proposition \ref{prop:bij-set} to $(E, \nu_l, F^{(2)})$, $\mathcal{A}$ and $A$, identified with $A \times A$.
\end{proof}

\begin{remark}\label{rem:triv}
    Another corollary of Proposition \ref{prop:bij-set} is that every ergodic type III CBER is isomorphic to its restriction on every positive measure subset $A \subset X$. Here, an isomorphism between two mcp CBERs $(X, \mu, E)$ and $(Y, \nu, F)$ is a Borel bijection $f:X \rightarrow Y$ preserving the measure class for which
    \[ x \mathbin{E} y \longleftrightarrow f(x) \mathbin{F} f(y).\]
    
    In particular, any isomorphism invariant $\Omega(\cdot) \in [0,\infty]$ of mcp CBERs that scales naturally, i.e,
    \[ \Omega(E, \mu) = \mu(A) \cdot \Omega\left(E|_A, \frac{\mu|_A}{\mu(A)}\right)\]
    (similarly to $\mathrm{cost}(E) - 1$ or the first $\ell^2$-Betti number in the pmp context), will be either $0$ or $\infty$. 
\end{remark}

\begin{theorem}\label{thm:graphslid}
    Let $(X, \mu, G)$ be an ergodic type III directed Borel graph and $H < G$ be an ergodic type III Borel strict subgraph. Let $g_{i}, g_{o} : X \rightarrow \N \cup \{\infty\}$ be Borel maps such that for all $x\in X$:
    \begin{align*}
        d_i(x) &\geq \ind_{H}(x) \\
        d_o(x) &\geq \outd_{H}(x)
    \end{align*}
    and each inequality is strict on a positive measure set. Then there exists a Borel graph $G' > H$ satisfying
    \begin{itemize}
        \item $G'$ has the same connected components as $G$:  $E_G = E_{G'}$ a.e.
        \item The in-degree of $G'$ is exactly $d_i$ : $\ind_{G'} \equiv d_i$.
        \item The out-degree of $G'$ is exactly $d_o$ : $\outd_{ G'} \equiv d_o$.
    \end{itemize}
    Furthermore, if $G$ is acyclic, then $G'$ can be chosen to be acyclic.
\end{theorem}
\begin{proof}
    Without loss of generality, suppose that for all $(x,y) \in G - H$, then $(x,y) \not\in E_H$ and there are no $x' \in [x]_H$ or $y' \in [y]_H$ with $(x', y), (y',x) \in G - H$. This is done by removing edges from $G$, which can be done in a Borel way using the Luzin-Novikov theorem.

    We first start by sliding the endpoints of edges to ensure the in-degree condition is satisfied. Consider the essential supremum $S = \sup_\mu \left(d_i - \ind_{H} \right) \in \{1, 2, ..., \infty \}$. First, let $\overrightarrow{G - H}$ be an orientation of $G-H$ as in Lemma \ref{lem:order}. Take a partition
    \begin{align*}  
        \overrightarrow{G - H} = \mathcal{A}_1 \sqcup \dots \sqcup \mathcal{A}_{S} \hspace{4em}&\text{ if } S \in \N,\\
        \overrightarrow{G - H} = \mathcal{A}_1 \sqcup \dots \sqcup \mathcal{A}_k \sqcup \dots \hspace{2em}&\text{ if } S = \infty,
        \end{align*}
    
    into non-null sets of oriented edges. Apply Lemma \ref{lem:bij-edge} repeatedly with $\mathcal{A}_k$ and 
    \[A_k = \left\{x \in X: d_i(x) \geq \ind_{H}(x) + k\right \}\] to find a Borel bijection $\phi_k: \mathcal{A}_k \rightarrow A_k$ such that $\phi_k(e) \mathbin{E}_H s(e)$ for all $e \in \mathcal{A}$. Let
    \[G'' = H\sqcup  \bigsqcup_{k = 1}^S\left\{ (\phi_k(e), t(e))  :  e \in \mathcal{A}_k \right\}\]
    The union is disjoint by our hypothesis, since there can be at most one edge of $G - H$ between the second coordinate $t(e)$ and the $H$-class $[s(e)]_H$ of the first coordinate. This graph still has the same connected components as $G$.
    
    We have that 
    \begin{align*}
    \ind_{G''}(x) &= \ind_{H}(x) + \# \{k: x \in A_k \} \\ 
    &= \ind_{H}(x) + \max\{k: x \in A_k \}\\
    &= \ind_{H}(x) + d_i(x) - \ind_{H}(x)\\
    &= d_i(x).
    \end{align*}
  		Reversing the orientation $\overrightarrow{G'' - H}$ and replacing $d_i$ with $d_o$, we repeat the argument and also obtain the desired equality for outdegrees.

    For the last point on acyclic graphs, notice given that $G$ is acyclic, the quotient tree $G/H$ is left unchanged by the edge slides above, and the number of edges between two $H$-components ($0$ or $1$ in the acyclic case) does not change either.
\end{proof}

\begin{remark}\label{rem:ergodic}
Ergodicity of $G$ and $H$ is not required, nor is $H$ required to be relatively ergodic in $G$. One could replace them by asking that both inequalities be strict on $H$-complete sections instead.
\end{remark}

\section{Antirigidity of trees in type III combinatorics}

\begin{theorem}\label{theo:main}
    Let $(X, \mu, G)$ be an ergodic type III non-amenable locally countable Borel acyclic graph and $d \in \{2, \dots, \infty\}$. Then there is a free action of $\F_d$ on $X$ such that $E_G = E_{\F_d}$. Further, the Schreier graph of this action does not have vanishing ends.
\end{theorem}
\begin{proof}
        By \cite[Theorem 1.3]{Tserunyan_2022}, find an amenable ergodic subforest $H < G$. The inclusion is strict since $G$ is non-amenable. If $H$ is type III, we can skip further down the proof, hence suppose it is type II. Without loss of generality, assume that $\mu$ is preserved by $H$. 
        
        For $\alpha > 0$, consider the set 
        \[\mathcal{D}' = \left\{(y,x) \in G - H: \w(y,x) \in (\alpha, \infty)\right\}.\] Let $\epsilon$ be small enough so that this set is non-null. Using \cite[Lemma 7.3]{kechris_topics_2004} and proper ergodicity, take a non-null subset $\mathcal{D} \subset \mathcal{D}'$ such that: 
        \begin{itemize}
            \item for $e, e' \in \mathcal{D}$, their endpoints do not overlap. 
            \item $(G - H) - \mathcal{D}$ is non-null.
        \end{itemize}
        Thus, $H' = H \sqcup \mathcal{D} < G$ is a strict subforest of $G$. We now show that $H'$ is type III. 

        Let $A \subset X$ be a complete section and $x \in A$. Let $M > 0$ be arbitrarily large and let $n \geq \frac{\ln M}{\epsilon}$. Inductively find a path $(x= x_0, y_0, x_1, y_1, ..., x_{n-1}, y_{n-1}, x_n)$ in $H'$ with $x_i \mathbin{E}_H y_i$ and $(y_{i -1}, x_i) \in \mathcal{D}$ with $\w( x_i, y_{i -1 }) > e^{\epsilon}$. This path always exists by ergodicity of $H$, ensuring that there is always a lighter end of an edge of $\mathcal{D}$ in $[x_i]_H$. Indeed, the path of lighter endpoint of edges in $\mathcal{D}$ being of positive measure, it must intersect every $H$-class, since $H$ is ergodic.
        
        Then,
        \begin{align*}
        \w( x_n, x_0)&= \w(y_0, x_0) \cdot \w( x_1, y_0)  \dots  \w(x_n,y_{n-1}) \\ 
        &> 1 \cdot e^\epsilon \dots  e^\epsilon  = e^{n\epsilon} > M.
        \end{align*}
        Since $H$ is ergodic, there exists $x' \in [x_n]_H \cap A$. Therefore, $\w$ is unbounded above on $A$, as $\w(x',x) > M \cdot \w( x', x_n) = M$. A similar argument shows unboundedness below. Thus, by Corollary \ref{cor:unbound}, $H'$ is type III.

        Without losing generality (since we only care about $E_G$), we can assume $H$ is generated by a single bijection $T_H$ on $X$. We then have a type III subforest $H' < G$ with the property that $\deg_{H'}(x) \leq 3$ for all $x \in X$.

        The first generator of the $\F_d$ action will be given by $T_H$. For the other generators, orient arbitrarily edges of $G - H - \mathcal{D}$, using Lemma \ref{lem:order}. Partition into non-null sets
        \begin{align*}  
        \overrightarrow{G - H} = A_1 \sqcup \dots \sqcup A_{d-1} \hspace{3.3em}&\text{ if } d \in \N,\\
        \overrightarrow{G - H} = A_1 \sqcup \dots \sqcup A_i \sqcup \dots \hspace{2em}&\text{ if } d = \infty,
        \end{align*}
        with $\mathcal{D} \subset A_1$. The edges of $A_i$ will be inductively slid to the generator $i + 1$. Apply Theorem \ref{thm:graphslid} using $H' < H' \sqcup A_1$ and $d_i, d_o \equiv 2$ to get a subtree $H' < H^{(1)}$ which is induced by a free $\F_2$-action, as it is a disjoint union of $H$ and a $\Z$-line, which can be taken without loss of generality to be generated by a $\Z$-action.
        
        For $d \geq 3$, one repeats this procedure with $H^{(1)} < H^{(1)} \sqcup A_2$ to get an $\F_3$ action, and so on and so forth.

        To see that the graph does not have vanishing ends, notice that every set of the form $N_{y,x}$ contains an $H'$-class. Since $H'$ is type III, the cocycle is then unbounded in $N_{y,x}$.
\end{proof}

\begin{remark}
As noted in Remark \ref{rem:ergodic}, $H$ needs not be ergodic in order to apply Theorem \ref{thm:graphslid}. However, it is unclear how to construct such an $H$ which would not be relatively ergodic in $G$, but for which both the set of lighter (resp. heavier) end of edges of $\mathcal{D}'$ are $H$-complete sections. In this sense, we use the full strength of \cite{Tserunyan_2022}.
\end{remark}

We now give a sketch of how to construct an example showing how type III is required, as opposed to compressibility, for Theorem 4.1.

\begin{example}\label{examp:IIinf}
    Consider a free pmp action $\F_2 \act (X, \mu)$ and the relation on $X \times \N$ given by
    \[ (x,n) \mathrel{E} (y,n) \Longleftrightarrow y \in \F_2 \cdot x \]
    This is a type II relation preserving an infinite measure. Choose $\nu$ a probability measure equivalent to the product of $\mu$ and the counting measure. Consider an action of $\F_\infty$ generating $E$ and let $Y \subset X \times \N$ be the convex hull of $X \times \{ 0 \}$ in the Schreier graph of the action. The subgraph induced on $Y$ must be locally finite, else one could find a treeing of the action $\F_2 \act X$ with infinite cost. In particular, the retract along the Schreier graph $X \times \N \rightarrow Y$ has fibers over $y$ a locally infinite tree almost surely. In particular, the subgraph induced by the retract is smooth, hence has only vanishing ends in $\nu$, as is proved in \cite{TTD}.
\end{example}

\bibliography{references}
\bibliographystyle{alpha}

\end{document}